\input amstex
\magnification=\magstep1 
\baselineskip=13pt
\documentstyle{amsppt}
\vsize=8.7truein \CenteredTagsOnSplits \NoRunningHeads
\def\UU{\Cal U}
 \def\den{\operatorname{Density}}
 \topmatter
 
\title Computing the partition function for cliques in a graph  \endtitle 
\author Alexander Barvinok \endauthor
\address Department of Mathematics, University of Michigan, Ann Arbor,
MI 48109-1043, USA \endaddress
\email barvinok$\@$umich.edu \endemail
\date October 2014 \enddate
\thanks  This research was partially supported by NSF Grants DMS 0856640 and DMS 1361541.
\endthanks 
\keywords algorithm, clique, partition function, subgraph density, graph  \endkeywords
\abstract We present a deterministic algorithm which, given a graph $G$ with $n$ vertices and an integer $1<m \leq n$, 
computes in $n^{O(\ln m)}$ time the sum of weights $w(S)$ over all $m$-subsets $S$ of the set of vertices of $G$, where 
$w(S)=\exp\left\{\gamma t m +O\left({1 \over m}\right)\right\}$ provided exactly $t{m \choose 2}$ pairs of vertices of $S$ span an edge of $G$ for some $0 \leq t \leq 1$.
Here $\gamma >0$ is an absolute constant: we can choose $\gamma=0.06$, and if
$n \geq 4m$ and $m \geq 10$, we can choose $\gamma=0.18$. This allows us to tell apart the graphs that do not have $m$-subsets of high density from the graphs that have sufficiently many $m$-subsets of high density, even when the probability to hit such a subset at random is exponentially small in $m$.
 \endabstract
\subjclass 15A15, 68C25, 68W25, 60C05  \endsubjclass

\endtopmatter

\document

\head 1. Introduction and main results \endhead

\subhead (1.1) Density of a subset in a graph \endsubhead
Let $G=(V, E)$ be an undirected graph with set $V$ of vertices and set $E$ of edges, without loops or multiple edges.
Let $S \subset V$ be a subset of the set of vertices of $G$. We define the {\it density} $\sigma(S)$ of $S$ as the ratio of the number of the edges of $G$ with both endpoints in $S$ to the maximum possible number of edges spanned by vertices in $S$:
$$\sigma(S)={\left| \{u, v\} \in E: \quad u, v \in S \right| \over {|S| \choose 2}}.$$
Hence 
$$0 \ \leq \ \sigma(S) \ \leq 1 \quad \text{for all} \quad S \subset V.$$ 
In particular, $\sigma(S)=0$ if and only if $S$ is an {\it independent set}, that is, no two vertices of $S$ span an edge of $G$ and $\sigma(S)=1$ if and only if $S$ is a {\it clique}, that is, every two vertices of $S$ span an edge of $G$. 

In this paper, we suggest a new approach to testing the existence of subsets $S$ of a given size $m=|S|$ with high density $\sigma(S)$. Namely, we present a deterministic algorithm, which, given a graph $G$ with $n$ vertices and an integer $1 < m \leq n$ computes in $n^{O(\ln m)}$ time (within a relative error of $0.1$, say) the sum 
$$\den_m(G)=\sum \Sb S \subset V \\ |S|=m \endSb \exp\left\{ \gamma m \sigma(S) + O\left({1 \over m} \right)\right\},\tag1.1.1$$
where  $\gamma>0$ is an absolute constant: we can choose $\gamma=0.06$ and if $n \geq 4m$ and $m \geq 10$ 
then we can choose $\gamma=0.18$. All implicit constants in the ``$O$" notation are absolute.

Let us fix two non-negative real numbers $\sigma$ and $\epsilon$ such that $\sigma + \epsilon \leq 1$. If there are no $m$-subsets $S$ with density $\sigma$ or higher then 
$$\den_m(G) \ \leq \ {n \choose m} \exp\left\{ \gamma m \sigma +O\left({1 \over m}\right)\right\}.$$
If, however, the there are sufficiently many $m$-subsets $S$ with density $\sigma + \epsilon$ or higher, that is, if the probability to hit such a subset at random is at least $2 \exp\{-\gamma \epsilon m\}$, then 
$$\den_m(G) \ \geq \ 2{n \choose m} \exp\left\{\gamma m \sigma  + O\left({1 \over m} \right)\right\}$$
and we can distinguish between these two cases in $n^{O(\ln m)}$ time. Note that ``many subsets" still allows for the probability to hit such a subset at random to be exponentially small in $m$. 

It turns out that we can compute in $n^{O(\ln m)}$ time an $m$-subset $S$, which is almost as dense as the average under the exponential weighting of (1.1.1), that is, an $m$-subset $S$ satisfying 
$$\exp\bigl\{\gamma m \sigma(S) \bigr\} \ \geq \ {1 \over 2} {n \choose m}^{-1} \den_m(G),  \tag1.1.2$$
cf. Remark 3.6.

We compute (1.1.1) using partition functions.

\subhead (1.2) The partition function of cliques \endsubhead Let $W=\left(w_{ij}\right)$ be an $n \times n$ symmetric matrix, interpreted as a matrix of weights on the edges of the complete graph on $n$ vertices. For an integer $1 < m \leq n$,
we define a polynomial
$$P_m(W)=\sum \Sb S \subset \{1, \ldots, n\} \\ |S|=m \endSb \prod \Sb \{i, j\} \subset S \\ i \ne j \endSb w_{ij},$$
which we call the {\it partition function} of cliques (note that it is different from what is known as the partition function of independent sets, see
\cite{SS05} and \cite{BG08}).
Thus if $W$ is the adjacency matrix of a graph $G$ with set $V=\{1, \ldots, n\}$ of vertices, so that
$$w_{ij}=\cases 1 &\text{if\ } \{i, j\} \in E \\ 0 &\text{otherwise}, \endcases$$
then $P_m(W)$ is the number of cliques in $G$ of size $m$.
Hence computing $P_m(W)$ is at least as hard as counting cliques. Let us choose an $0 < \alpha < 1$ and modify the weight $W$ 
as follows:
$$w_{ij}=\cases 1 &\text{if\ } \{i, j\} \in E \\ \alpha&\text{otherwise.} \endcases$$
Then $P_m(W)$ counts all cliques of size $m$ in the complete graph on $n$ vertices, only that a clique containing exactly $k$ non-edges of $G$ is weighted down by a factor of $\alpha^k$. In other words, an $m$-subset $S \subset V$ is counted with weight 
$$\exp\left\{ \bigl(1-\sigma(S)\bigr) {m \choose 2} \ln \alpha \right\}.$$

In this paper, we present a deterministic algorithm, which, given an $n \times n$ symmetric matrix $W=\left(w_{ij}\right)$ such that 
$$\left|w_{ij}-1 \right| \ \leq \ {\gamma \over m-1} \quad \text{for all} \quad i,j \tag1.2.1$$
computes the value of $P_m(W)$ within relative error $\epsilon$ in $n^{O\left(\ln m - \ln \epsilon\right)}$ time. Here 
$\gamma >0$ is an absolute constant; we can choose $\gamma=0.06$ and if $n \geq 4m$ and $m \geq 10$, we can choose 
$\gamma =0.18$.

\subhead (1.3) The idea of the algorithm \endsubhead
Let $J$ denote the $n \times n$ matrix filled with 1s. Given an $n \times n$ symmetric matrix $W$, we consider the univariate function 
$$f(t)=\ln P_m \bigl(J + t(W-J)\bigr). \tag1.3.1$$
Clearly,
$$f(0)=\ln P_m(J) =\ln {n \choose m} \quad \text{and} \quad f(1)=\ln P_m(W).$$
Hence our goal is to approximate $f(1)$ and we do it by using the Taylor polynomial expansion of $f$ at $t=0$:
$$f(1) \approx f(0) +\sum_{k=1}^l {1 \over k!} {d^k \over dt^k} f(t) \Big|_{t=0}. \tag1.3.2$$
It turns out that the right hand side of (1.3.2) can be computed in $n^{O(l)}$ time.
We present the algorithm in Section 2. The quality of the approximation (1.3.2) depends on the location of {\it complex} zeros of $P_m$.

\proclaim{(1.4) Lemma} Suppose that there exists a real $\beta>1$ such that 
$$P_m \bigl(J +z(W-J)\bigr) \ne 0 \quad \text{for all} \quad z \in {\Bbb C} \quad \text{satsifying} \quad |z| \leq \beta.$$
Then the right hand side of (1.3.2) approximates $f(1)$ within an additive error of 
$${m(m-1) \over 2(l+1)\beta^l(\beta-1)}.$$
\endproclaim
In particular, for a fixed $\beta >1$, to ensure an additive error of $0< \epsilon < 1$, we can choose 
$l=O\left(\ln m -\ln \epsilon\right)$, which results in the algorithm for approximating $P_m(W)$ within relative error $\epsilon$ in
$n^{O\left(\ln m -\ln \epsilon\right)}$ time. We prove Lemma 1.4 in Section 2.

Thus we have to identify a class of matrices $W$ for which the number $\beta >1$ of Lemma 1.4 exists.
We prove the following result.

\proclaim{(1.5) Theorem} There is an absolute constant $\omega >0$  (one can choose $\omega=0.061$, and if 
$m \geq 10$ and $n \geq 4m$, one can choose $\omega=0.181$) such that for any positive integers $1 < m \leq n$, for any $n \times n$ symmetric complex matrix $Z=\left(z_{ij}\right)$, we have 
$$P_m(Z) \ne 0$$
as long as
$$\left| z_{ij}-1\right| \ \leq \ {\omega \over m-1} \quad \text{for all} \quad i,j.$$
\endproclaim

We prove Theorem 1.5 in Section 3. 
Theorem 1.5 implies that if (1.2.1) holds, 
where $0 < \gamma < \omega$ is an absolute constant, we can choose $\beta=\omega/\gamma$ in Lemma 1.4
 and obtain an algorithm which computes $P_m(W)$ within a given relative error $\epsilon$ in $n^{O(\ln m-\ln \epsilon)}$ time. 
 Thus we can choose $\gamma=0.06$ and $\beta=61/60$ and if $m \geq 10$ and $n \geq 4m$, we can choose 
 $\gamma=0.18$ and $\beta=181/180$.

\subhead (1.6) Weighted enumeration of subsets \endsubhead
Given a graph $G$ with set of vertices $V=\{1, \ldots, n\}$ and set $E$ of edges, we define a weight $W$ by 
$$w_{ij}=\cases 1+{\gamma \over m-1} &\text{if \ } \{i, j \} \in E \\ 1-{\gamma \over m-1} &\text{otherwise,} \endcases$$
where $\gamma>0$ is an absolute constant in (1.2.1). Thus we can choose $\gamma=0.06$ and if $m \geq 10$ and
$n \geq 4m$, we can choose $\gamma=0.18$. We define the functional (1.1.1) by 
$$\den_m(G)= e^{\gamma m} \left( 1+{\gamma \over m-1}\right)^{-{m \choose 2}}P_m(W).$$
Then we have 
$$\operatorname\den_m(G)=\sum \Sb S \subset V \\ |S|=m \endSb w(S),$$
where 
$$\split w(S)=&e^{\gamma m} \left(1+{\gamma \over m-1}\right)^{(t-1){m \choose 2}}\left(1-{\gamma \over m-1}\right)^{(1-t){m \choose 2}} \\= 
&\exp\left\{\gamma m t + O\left({1 \over m}\right) \right\}, \endsplit$$
provided exactly $t{m \choose 2}$ pairs of vertices of $S$ span an edge of $G$, that is, provided 
$\sigma(S)=t$.
Thus we are in the situation described in Section 1.1.

\subhead (1.7) Comparison with results in the literature \endsubhead
The topic of this paper touches upon some very active research areas, such as finding dense subgraphs in a given graph (see, for example,  \cite{Bh12} and references therein), computational complexity of partition functions (see, for example, \cite{BG05} and references therein) and complex zeros of partition functions (see, for example, \cite{SS05} and references therein). Detecting dense subsets is notoriously hard. It is an NP-hard problem to approximate the size $|S|$ of the largest clique $S$ within a factor of $|V|^{1-\epsilon}$ for any fixed $0 < \epsilon \leq 1$ \cite{H\aa99}, \cite{Zu07}.
In \cite{F+01}, a polynomial time algorithm for approximating the highest density of an $m$-subset within a factor of 
$O(n^{{1 \over 3} -\epsilon})$ for some small $\epsilon>0$ was constructed. The paper \cite{B+10} presents an algorithm of $n^{O(1/\epsilon)}$ complexity which approximates the highest density of an $m$-subset within a factor of $n^{{1 \over 4} +\epsilon}$ and an algorithm of $O\left(n^{\ln n}\right)$ complexity which approximates the density within a factor of $O(n^{1/4})$, the current record (note that \cite{F+01}, \cite{Bh12} and \cite{B+10} normalize density slightly differently than we do). It has also been established that modulo some plausible complexity assumptions, it is hard to approximate the highest density of an $m$-subset within a constant factor, fixed in advance, see \cite{Bh12}.

We note that while searching for a dense $m$-subset in a graph on $n$ vertices, the most interesting case to consider is that of a growing $m$ such that $m=o(n)$, since for any fixed $\epsilon>0$ one can compute in polynomial time the maximum number of edges of $G$ spanned by an $m$-subset $S$ of vertices of $G$ within an additive error of $\epsilon n^2$ \cite{FK99} (the complexity of the algorithm is exponential in $1/\epsilon$). Of course, for any constant $m$, fixed in advance, all subsets of size $m$ can be found by the exhaustive search in polynomial time.

The appearance of partition functions in combinatorics can be traced back to the work of Kasteleyn on the statistics of dimers (perfect matchings), see for example, Section 8.7 of \cite{LP09}. A randomized polynomial time approximation algorithm (based on the Markov Chain Monte Carlo approach) for computing the partition function of the classical Ising model was constructed by Jerrum and Sinclair in \cite{JS93}.
A lot of work has recently been done on understanding the computational complexity of the partition function of graph homomorphisms \cite{BG05}, where a certain ``dichotomy" of instances has been established: it is either $\#P$-hard to compute exactly (in most interesting cases) or it is computable in polynomial time exactly (in few cases). Recently, an approach inspired by the ``correlation decay" phenomenon in statistical physics was used to construct deterministic approximation algorithms for partition functions in combinatorial problems \cite{BG08}, \cite{B+07}. To the best of author's knowledge, the partition function $P_m(W)$ introduced in Section 1.2 of this paper has not been studied before and the idea of our algorithm sketched in Section 1.3 is also new. On the other hand, one can view our method as inspired by the intuition from statistical physics. 
Essentially, the algorithm computes the partition function similar to the one used in the ``simulated annealing" method, see \cite{Az88}, for the temperature that is provably higher than the phase transition threshold (the role of the ``temperature" is played by $1/\gamma m$, where $\gamma$ is the constant in (1.1.1)). The algorithm uses the fact that the partition function has nice analytic properties at temperatures above the phase transition temperature. 

The corollary asserting that graphs without dense $m$-subsets can be efficiently separated from graphs with many dense $m$-subsets is vaguely similar in spirit to the property testing approach \cite{G+98}.

The result requiring most work is Theorem 1.5 bounding the complex roots of the partition function away from the matrix of all 1s. 
Studying complex roots of partition functions was initiated by Lee and Yang \cite{LY52} and it is a classical topic by now, because of its importance to locating the phase transition, see, for example, \cite{SS05} and Section 8.5 of \cite{LP09}.

\subhead (1.8) Open question \endsubhead Is it true that for {\it any} $\gamma >0$, fixed in advance, the density functional (1.1.1) can be computed (within a relative error of $0.1$, say) in $n^{O(\ln m)}$ time? 

\head 2. The algorithm \endhead

\subhead (2.1) The algorithm for approximating the partition function \endsubhead Given an $n \times n$ symmetric matrix
$W=\left(w_{ij}\right)$, we present an algorithm that computes the right hand side of (1.3.2) 
for the function $f(t)$ defined by (1.3.1).

Let 
$$g(t)=P_m \bigl(J+t(W-J)\bigr) =
\sum \Sb S \subset \{1, \ldots, n\} \\ |S|=m \endSb \prod \Sb \{i, j\} \subset S \\ i \ne j \endSb \bigl(1+t(w_{ij}-1)\bigr)  , \tag2.1.1$$
so $f(t)=\ln g(t)$.
Hence
$$f'(t)={g'(t) \over g(t)} \quad \text{and} \quad g'(t)=g(t) f'(t).$$
Therefore, for $k \geq 1$ we have 
$${d^k \over dt^k}g(t)\Big|_{t=0} =\sum_{j=0}^{k-1} {k-1 \choose j} \left( {d^j \over dt^j} g(t) \Big|_{t=0} \right) 
\left({d^{k-j} \over dt^{k-j}} f (t) \Big|_{t=0}\right) \tag2.1.2$$
(we agree that the 0-th derivative of $g$ is $g$).

We note that $g(0)={n \choose m}$. If we compute the values of 
$${d^k \over dt^k} g(t)\Big|_{t=0} \quad \text{for} \quad k=1, \ldots, l, \tag2.1.3$$
then the formulas (2.1.2) for $k=1, \ldots, l$ provide a non-degenerate triangular system of linear equations that allows us to 
compute 
$${d^k \over dt^k} f(t)\Big|_{t=0} \quad \text{for} \quad k=1, \ldots, l.$$
Hence our goal is to compute the values (2.1.3).

We have 
$$\split {d^k \over dt^k} g(t) \Big|_{t=0} =&\sum \Sb S \subset \{1, \ldots, n\} \\ |S|=m \endSb \sum 
\Sb I=\left(\{i_1, j_1\}, \ldots, \{i_k, j_k\}\right) \\
\{i_1, j_1\}, \ldots, \{i_k, j_k\}\subset S \endSb \left(w_{i_1 j_1}-1\right) \cdots \left(w_{i_k j_k} -1 \right), \endsplit$$
where the inner sum is taken over all ordered sets $I$ of $k$ pairs $\{i_1, j_1\}$, $\ldots$, $\{i_k, j_k\}$ of points from $S$.

For an ordered set
$$I=\left( \{i_1, j_1\}, \ldots, \{i_k, j_k\}\right)$$ 
of $k$ pairs of points from the set $\{1, \ldots, n\}$, let
$$\rho(I)=\left| \bigcup_{s=1}^k \{i_s, j_s\} \right|$$
be the total number of distinct points among $i_1, j_1, \ldots, i_k, j_k$.
Since there are exactly 
$${n-\rho(I) \choose m - \rho(I)}$$
$m$-subsets $S  \subset \{1, \ldots, n\}$ containing the pairs $\{i_1, j_1\}$, $\ldots$, $\{i_k, j_k\}$,
we can further rewrite
$$\split &{d^k \over dt^k} g(t)\Big|_{t=0}\\ &\quad=\sum \Sb  I=\left(\{i_1, j_1\}, \ldots, \{i_k, j_k\} \right)
\\ \{i_1, j_1\}, \ldots, \{i_k, j_k\}  \subset \{1, \ldots, n\} \endSb 
{n-\rho(I) \choose m - \rho(I)}
\left(w_{i_1 j_1}-1\right) \cdots \left(w_{i_k j_k} -1 \right), \endsplit$$
where the sum is taken all ordered sets $I$ of $k$ pairs $\{i_1, j_1\}, \ldots, \{i_k, j_k\}$ of points from the set 
$\{1, \ldots, n\}$. 
The algorithm consists in computing the latter sum. Since the sum contains not more than 
$n^{2k}=n^{O(l)}$ terms, the complexity of the algorithm is indeed $n^{O(l)}$, as claimed.

\subhead (2.2) Proof of Lemma 1.4 \endsubhead The function $g(z)$ defined by (2.1.1) is a polynomial in $z$ of degree $d \leq {m \choose 2}$ with $g(0)={n \choose m} \ne 0$, so we factor
$$g(z)= g(0) \prod_{i=1}^d \left(1-{z \over \alpha_i}\right)
,$$
where $\alpha_1, \ldots, \alpha_d$ are the roots of $g(z)$. By the condition of Lemma 1.4, 
we have 
$$\left|\alpha_i\right| \ \geq \ \beta >1 \quad \text{for} \quad i=1, \ldots, d.$$
Therefore,
$$f(z)=\ln g(z)= \ln g(0) + \sum_{i=1}^d \ln \left(1-{z \over \alpha_i}\right) \quad \text{provided} \quad |z| \leq 1 , \tag2.2.1$$
where we choose the branch of $\ln g(z)$ that is real at $z=0$. Using the standard Taylor expansion, we obtain 
$$\ln \left(1-{1 \over \alpha_i}\right) =-\sum_{k=1}^l {1 \over k} \left({1 \over \alpha_i}\right)^k + \zeta_l,$$
where
$$\left|\zeta_l\right|= \left| \sum_{k=l+1}^{+\infty} {1 \over k} \left( {1 \over \alpha_i} \right)^k \right| \ \leq \ 
{1 \over (l+1) \beta^l (\beta-1)}.$$
Therefore, from (2.2.1) we obtain
$$f(1)=f(0)  + \sum_{k=1}^l \left(-{1 \over k} \sum_{i=1}^d \left({1 \over \alpha_i}\right)^k\right)  +\eta_l,$$
where 
$$\left| \eta_l \right| \ \leq \ {m(m-1)\over 2(l+1)\beta^l (\beta-1)}.$$
It remains to notice that 
$$-{1 \over k} \sum_{i=1}^d \left({1 \over \alpha_i}\right)^k ={1 \over k!} {d^k \over dt^k} f(t) \Big|_{t=0}.$$
{\hfill \hfill \hfill} \qed

\head 3. Proof of Theorem 1.5 \endhead

\definition{(3.1) Definitions} For $\delta>0$, we denote  by $\UU(\delta) \subset {\Bbb C}^{n(n+1)/2}$ the closed polydisc 
$$\UU(\delta) =\Bigl\{ Z=\left(z_{ij} \right): \quad z_{ij}=z_{ji} \quad \text{and} \quad  \left| z_{ij}-1 \right| \leq \delta \quad \text{for all} \quad i, j \Bigr\}.$$
For a subset $\Omega \subset \{1, \ldots, n\}$ where $0 \leq |\Omega| \leq m$, we define
$$P_{\Omega}(Z)=\sum \Sb S \subset \{1, \ldots, n\} \\ |S|=m \\ S \supset \Omega \endSb 
\prod \Sb \{i, j\} \subset S \\ i \ne j \endSb z_{ij}. $$
In words: $P_{\Omega}(Z)$ is the restriction of the sum defining the clique partition function to the subsets $S$ that contain a given set 
$\Omega$.
In particular,
$$P_{\emptyset}(Z)=P_m(Z),$$
and if $|\Omega| < m$ then 
$$P_{\Omega}(Z) ={1 \over m-|\Omega|}  \sum_{i :\ i \notin \Omega} P_{\Omega \cup \{i\}} (Z). \tag3.1.1$$
We note that the natural action of the symmetric group $S_n$ on $\{1, \ldots, n\}$ induces the natural action of $S_n$ on the space of polynomials in $Z$ which permutes the polynomials $P_{\Omega}(Z)$.
\enddefinition

First, we establish a simple geometric lemma.

\proclaim{(3.2) Lemma} Let $u_1, \ldots, u_n \in {\Bbb R}^d$ be non-zero vectors such that for some $0 \leq \alpha < \pi/2$
the angle between any two vectors $u_i$ and $u_j$ does not exceed $\alpha$. Let 
$u=u_1 + \ldots + u_n$. Then 
$$\|u\| \ \geq \ \sqrt{\cos \alpha}  \sum_{i=1}^n \|u_i\|.$$
\endproclaim
\demo{Proof} We have 
$$\|u\|^2 =\sum_{1 \leq i, j \leq n} \langle u_i, u_j \rangle \ \geq \ \sum_{1 \leq i, j \leq n} \|u_i\| \|u_j\| \cos \alpha =
\left(\cos \alpha \right) \left( \sum_{i=1}^n \|u_i \| \right)^2.$$
{\hfill \hfill \hfill} \qed
\enddemo

We prove Theorem 1.5 by the reverse induction on $|\Omega|$, using Lemma 3.2 and the following two lemmas.

\proclaim{(3.3) Lemma} Let us fix real $0 < \tau < 1$, real $\delta > 0$ and an integer $1 \leq  r \leq m$.
Suppose that for all $\Omega \subset \{1, \ldots, n \}$ such that $|\Omega|=r$, 
for all $Z \in \UU(\delta)$, we have 
$P_{\Omega}(Z) \ne 0$ 
and that for all $i=1, \ldots, n$ we have
$$\left|P_{\Omega}(Z)\right| \ \geq \ {\tau \over m-1} \sum_{j=1}^n \left|z_{ij}\right| \left| {\partial \over \partial z_{ij}} P_{\Omega}(Z) \right|.$$
Then, for any pair of subsets 
$$\Omega_1, \Omega_2 \subset \{1, \ldots, n\} \quad \text{such that} \quad \left|\Omega_1\right|=\left|\Omega_2\right|=r \quad \text{and} \quad \left| \Omega_1 \Delta \Omega_2 \right|=2,$$
and for any $Z \in \UU(\delta)$, the angle between the complex numbers $P_{\Omega_1}(Z)$ and $P_{\Omega_2}(Z)$,
interpreted as vectors in ${\Bbb R}^2 ={\Bbb C}$, does not exceed
$$\theta={ 4\delta (m-1) \over \tau(1-\delta)}$$
and the ratio of $\left|P_{\Omega_1}(Z)\right|$ and $\left|P_{\Omega_2}(Z)\right|$ does not exceed
$$\lambda =\exp\left\{{4 \delta (m-1) \over \tau(1-\delta)}\right\}.$$
\endproclaim
\demo{Proof} Let us choose an $\Omega \in \{1, \ldots, n\}$ such that $|\Omega|=r$.
Since $P_{\Omega}(Z) \ne 0$ for all $Z \in \UU(\delta)$, we can choose a branch of 
$\ln P_{\Omega}(Z)$ that is a real number when $Z=J$.
Then 
$${\partial \over \partial z_{ij}} \ln P_{\Omega}(Z)= \left({\partial \over \partial z_{ij}} P_{\Omega}(Z)\right) \Big/P_{\Omega}(Z).$$
Since $\left|z_{ij}\right| \geq 1-\delta$ for all $Z \in \UU(\delta)$, we have 
$$\sum_{j=1}^n \left| {\partial \over \partial z_{ij}} \ln P_{\Omega}(Z) \right| \ \leq \ {m-1 \over \tau(1-\delta)} \quad \text{for} \quad i=1, \ldots, n. \tag3.3.1$$
Without loss of generality, we assume that $1 \in \Omega_1$, $2 \notin \Omega_1$ and $\Omega_2=\Omega_1 \setminus \{1\}\cup\{2\}$. 

Given $A \in \UU(\delta)$, let $B \in \UU(\delta)$ be the matrix defined by 
$$b_{1j} =b_{j1}=a_{2j}=a_{j2}\quad \text{and} \quad b_{2j}=b_{j2}=a_{1j}=a_{j1} \quad \text{for} \quad j=1, \ldots, n$$
and 
$$b_{ij} =a_{ij} \quad \text{for all other} \quad i,j.$$
Then
$$P_{\Omega_2}(A)=P_{\Omega_1}(B)$$ and 
$$\split &\left| \ln P_{\Omega_1}(A) - \ln P_{\Omega_2}(A) \right|=
\left| \ln P_{\Omega_1}(A) - \ln P_{\Omega_1}(B) \right| \\ &\quad \leq \ \left( \sup_{Z \in \UU(\delta)} 
\sum_{j=1}^n \left| {\partial \over \partial z_{1j}} \ln P_{\Omega_1}(Z) \right| + \sum_{j=1}^n 
\left| {\partial \over \partial z_{2j}} \ln P_{\Omega_1}(Z) \right|  \right) \\
&\qquad \times \max_{j=1, \ldots, n} \left\{ \left| a_{1j} -b_{1j}\right|, \left|a_{2j}-b_{2j}\right|\right\}. \endsplit$$
Since $\left|a_{ij}-b_{ij}\right| \leq 2\delta$ for all $A, B \in \UU(\delta)$,
by (3.3.1) we conclude that the angle between the complex numbers $P_{\Omega_1}(A)$ and $P_{\Omega_2}(A)$ does not exceed $\theta$ and that the ratio of their absolute values does not exceed $\lambda$.
{\hfill \hfill \hfill} \qed
\enddemo

\proclaim{(3.4) Lemma} Let us fix real $\delta >0$, real $0 \leq \theta < \pi/2$, real $\lambda >1$ and an integer $1 < r \leq m$. 
Suppose that for all $\Omega \subset \{1, \ldots, n \}$ such that $m \geq |\Omega| \geq r$ and
for all $Z \in \UU(\delta)$, we have 
$P_{\Omega}(Z) \ne 0$  and that for any pair of subsets 
$$\Omega_1, \Omega_2 \subset \{1, \ldots, n\} \quad \text{such that} \quad m \ \geq \ \left|\Omega_1\right|=
\left|\Omega_2\right| \geq r \quad \text{and} \quad 
\left| \Omega_1 \Delta \Omega_2 \right|=2,$$
and for any $Z \in \UU(\delta)$,  the angle between two complex numbers $P_{\Omega_1}(Z)$ and $P_{\Omega_2}(Z)$ considered as vectors in ${\Bbb R}^2={\Bbb C}$ does not exceed $\theta$ while the ratio of $\left|P_{\Omega_1}(Z)\right|$ and $\left| P_{\Omega_2}(Z)\right|$ does not exceed $\lambda$. Let
$$\tau={\cos \theta \over \lambda}. $$
Then, for any subset $\Omega \subset \{1, \ldots, n\}$ such that $|\Omega|=r-1$, all  $Z \in \UU(\delta)$ and all
$i=1,\ldots, n$, we have 
$$\left|P_{\Omega}(Z)\right| \ \geq \ {\tau \over m-1} \sum_{j=1}^n \left|z_{ij}\right| \left| {\partial \over \partial z_{ij}} P_{\Omega}(Z)\right|. \tag3.4.1$$
In addition, 
$$\text{if} \quad {n \over m}  \ \geq \ {\lambda \over \sqrt{\cos \theta}},$$ the inequality (3.4.1) holds with 
$$\tau =\sqrt{\cos \theta}.$$
\endproclaim
\demo{Proof} Let us choose a subset $\Omega \subset \{1, \ldots, n\}$ such that $|\Omega|=r-1$. Let us define
$$\Omega_j=\Omega \cup \{j\} \quad \text{for} \quad j \notin \Omega.$$
Suppose first that $i \in \Omega$. Then, for $j \in \Omega$ we have 
$$ {\partial \over \partial z_{ij}} P_{\Omega}(Z) = z_{ij}^{-1} P_{\Omega}(Z)$$
and hence
$$\left| z_{ij}\right|  \left| {\partial \over \partial z_{ij}} P_{\Omega}(Z) \right| =\left| P_{\Omega}(Z) \right| \quad \text{provided}
\quad i,j \in \Omega, \quad i \ne j.$$
For $j \notin \Omega$, we have 
$${\partial \over \partial z_{ij}} P_{\Omega}(Z) ={\partial \over \partial z_{ij}} P_{\Omega_j}(Z) = 
z_{ij}^{-1} P_{\Omega_j}(Z)$$ 
and hence 
$$\left| z_{ij}\right|  \left| {\partial \over \partial z_{ij}} P_{\Omega}(Z) \right| =\left| P_{\Omega_j}(Z) \right| \quad \text{provided}
\quad i \in \Omega \quad \text{and} \quad j \notin \Omega.$$
Summarizing,
$$\sum_{j=1}^n \left|z_{ij}\right| \left| {\partial \over \partial z_{ij}} P_{\Omega}(Z)\right| =(r-2) \left|P_{\Omega}(Z)\right| +
\sum_{j \notin \Omega} \left|P_{\Omega_j}(Z)\right|. \tag3.4.2$$
On the other hand, by (3.1.1), we have 
$$P_{\Omega}(Z)={1 \over m-r+1} \sum_{j: \ j \notin \Omega} P_{\Omega_j}(Z).$$
By Lemma 3.2, we have 
$$\left| P_{\Omega}(Z)\right| \ \geq \ {\sqrt{\cos \theta} \over m-r+1} \sum_{j: \ j \notin \Omega} \left| P_{\Omega_j}(Z) \right|
\tag3.4.3$$
and
$$(m-1) \left|P_{\Omega}(Z)\right| \ \geq \ (r-2) \left|P_{\Omega}(Z)\right| +
\sqrt{\cos \theta} \sum_{j: \ j \notin \Omega} \left| P_{\Omega_j}(Z) \right|$$
and the inequality (3.4.1) with $\tau=\sqrt{\cos \theta}$ follows by (3.4.2).

Suppose now that $i \notin \Omega$. Then, for $j \in \Omega$, we have 
$${\partial \over \partial z_{ij}} P_{\Omega}(Z) ={\partial \over \partial z_{ij}} P_{\Omega_i}(Z) = 
z_{ij}^{-1} P_{\Omega_i}(Z)$$ 
and hence 
$$\left| z_{ij}\right|  \left| {\partial \over \partial z_{ij}} P_{\Omega}(Z) \right| =\left| P_{\Omega_i}(Z) \right| \quad \text{provided}
\quad i \notin \Omega \quad \text{and} \quad j \in \Omega. \tag3.4.4$$ 
If $r=m$ then 
$${\partial \over \partial z_{ij}} P_{\Omega}(Z) =0 \quad \text{for any} \quad j \notin \Omega $$
and
$$\sum_{j=1}^n \left|z_{ij}\right| \left| {\partial \over \partial z_{ij}} P_{\Omega}(Z)\right| =(m-1) \left|P_{\Omega_i}\right|
\quad \text{if} \quad r=m. $$
By (3.1.1), we have 
$$P_{\Omega}(Z)=\sum_{j \notin \Omega} P_{\Omega_j}(Z)$$ 
and hence by Lemma 3.2
$$\left| P_{\Omega}(Z)\right| \ \geq \ \sqrt{\cos \theta} \sum_{j \notin \Omega} \left| P_{\Omega_j}(Z) \right| \ \geq \ 
\sqrt{\cos \theta} \left|P_{\Omega_i}(Z)\right|,$$
which proves (3.4.1) in the case of $r=m$ with $\tau=\sqrt{\cos \theta}$.

If $r< m$, then for $j \notin \Omega_i$, let $\Omega_{ij}=\Omega \cup \{i, j\}$. Then, for $j \notin \Omega_i$, we have 
$${\partial \over \partial z_{ij}} P_{\Omega}(Z) ={\partial \over \partial z_{ij}} P_{\Omega_{ij}}(Z) = 
z_{ij}^{-1} P_{\Omega_{ij}}(Z)$$ 
and hence 
$$\left| z_{ij}\right|  \left| {\partial \over \partial z_{ij}} P_{\Omega}(Z) \right| =\left| P_{\Omega_{ij}}(Z) \right| \quad \text{provided}
\quad i \notin \Omega \quad \text{and} \quad j \notin \Omega.$$
From (3.4.4),
$$\sum_{j=1}^n \left|z_{ij}\right| \left| {\partial \over \partial z_{ij}} P_{\Omega}(Z)\right| =(r-1) \left|P_{\Omega_i}(Z)\right| +
\sum_{j \notin \Omega_i} \left|P_{\Omega_{ij}}(Z)\right| \quad \text{if} \quad r < m. \tag3.4.5$$
By (3.1.1), we have 
$$P_{\Omega_i}(Z)={1 \over m-r} \sum_{j:\ j \notin \Omega_i} P_{\Omega_{ij}}(Z)$$
and hence by Lemma 3.2, we have 
$$\left|P_{\Omega_i}(Z)\right| \ \geq \ {\sqrt{\cos \theta} \over m-r} \sum_{j:\ j \notin \Omega_i} \left|P_{\Omega_{ij}}(Z)\right|.$$
Comparing this with (3.4.5), we conclude as above that
$$\left|P_{\Omega_i}(Z)\right| \ \geq \ {\sqrt{\cos \theta} \over m-1} \sum_{j=1}^n \left| z_{ij}\right|  \left| {\partial \over \partial z_{ij}} P_{\Omega}(Z) \right|. \tag3.4.6$$
It remains to compare $\left|P_{\Omega_i}(Z)\right|$ and $\left|P_{\Omega}(Z)\right|$. 
Since the ratio of any two $\left|P_{\Omega_{j_1}}\right|$, $\left|P_{\Omega_{j_2}}\right|$ does not exceed $\lambda$, from
(3.4.3) we conclude that
$$\split \left|P_{\Omega}(Z)\right| \ \geq \ &{\sqrt{\cos \theta} \over \lambda}{n-r+1 \over m-r+1} \left|P_{\Omega_i}(Z)\right| 
\\ \geq \ &{\sqrt{\cos \theta} \over \lambda} \left|P_{\Omega_i}(Z)\right|\quad \text{for all}
\quad i \notin \Omega. \endsplit $$ 
The proof now follows by (3.4.6).
{\hfill \hfill \hfill} \qed
\enddemo

\subhead (3.5) Proof of Theorem 1.5 \endsubhead

One can see that for a sufficiently small $\omega >0$, the equation 
$$\theta ={4 \omega \exp\{\theta\} \over (1-\omega) \cos \theta}$$
has a solution $0 < \theta < \pi/2$. Numerical computations show that one can choose $\omega=0.061$, in which case
$$\theta \approx 0.4580097179.$$
Let 
$$\lambda=\exp\{\theta\} \approx 1.580924366 \quad \text{and} \quad \tau={\cos \theta \over \exp \{\theta\}} \approx 0.5673480171.$$
Let 
$$\delta={\omega \over m-1}.$$
For $r=m, \ldots,  1$ we prove the following three statements (3.5.1) -- (3.5.3):
\bigskip
(3.5.1) Let $\Omega \subset \{1, \ldots, n\}$ be a set such that $|\Omega|=r$. Then, for any $Z \in \UU(\delta)$, 
$Z=\left(z_{ij}\right)$, we have 
$$P_{\Omega}(Z) \ne 0;$$
\medskip
(3.5.2) Let $\Omega \subset \{1, \ldots, n\}$ be a set such that $|\Omega|=r$. Then for any $Z \in \UU(\delta)$ and
any $i=1, \ldots, n$, we have 
$$\left|P_{\Omega}(Z)\right| \ \geq \ {\tau \over m-1} \sum_{j=1}^n \left|z_{ij}\right| \left|{\partial \over \partial z_{ij}} 
P_{\Omega}(Z)\right|;$$
\medskip
(3.5.3) Let $\Omega_1, \Omega_2 \subset \{1, \ldots, n\}$ be sets such that 
$\left|\Omega_1\right| =\left| \Omega_2\right| =r$ and $\left|\Omega_1 \Delta \Omega_2\right|=2$.
Then the angle between the complex numbers $P_{\Omega_1}(Z)$ and $P_{\Omega_2}(Z)$, considered as vectors in 
${\Bbb R}^2={\Bbb C}$, does not exceed $\theta$, whereas the ratio of $\left|P_{\Omega_1}(Z)\right|$ and 
$\left|P_{\Omega_2}(Z)\right|$ does not exceed $\lambda$.
\bigskip
Suppose that $r=m$ and let $\Omega \subset \{1, \ldots, n\}$ be a subset such that $\left|\Omega\right|=m$.
Then 
$$P_{\Omega}(Z) = \prod \Sb \{i, j \} \subset \Omega \\ i \ne j \endSb z_{ij},$$
so clearly (3.5.1) holds for $r=m$.
Moreover, for $i \in \Omega$, we have
$$\sum_{j=1}^n \left|z_{ij}\right| \left| {\partial \over \partial z_{ij}} P_{\Omega}(Z) \right| = (m-1) \left| P_{\Omega}(Z)\right|,$$
while for $ i \notin \Omega$, we have 
$$\sum_{j=1}^n \left|z_{ij}\right| \left| {\partial \over \partial z_{ij}} P_{\Omega}(Z) \right| =0,$$
so (3.5.2) holds for $r=m$ as well. 

Lemma 3.3 implies that if (3.5.1) and (3.5.2) hold for sets $\Omega$ of cardinality $r$ then (3.5.3) holds of sets $\Omega$ of cardinality  $r$. Lemma 3.4 implies that if (3.5.1) and (3.5.3) hold for sets $\Omega$ of cardinality $r$ then (3.5.2) holds for sets $\Omega$ of cardinality $r-1$. Formula (3.1.1), Lemma 3.2 and statements (3.5.1) and (3.5.3) for sets $\Omega$ of cardinality $r$ imply statement (3.5.1) for sets $\Omega$ of cardinality $r-1$. Therefore, we conclude that the statements (3.5.1) and (3.5.3) hold for sets 
$\Omega$ such that $|\Omega|=1$. Formula (3.1.1) and Lemma 3.2 imply then that
$$P_{\emptyset}(Z)=P_m(Z) \ne 0,$$
as desired.

A more careful analysis establishes that if $m \geq 10$ then for $\omega=0.181$ the equation 
$$\theta ={4 \omega  \over \left(1-\dsize{\omega \over m-1}\right)\sqrt{ \cos \theta}}$$
has a solution 
$$0< \theta < 1.02831829$$
with 
$$\lambda = \exp\{\theta\} < 2.8 \quad \text{and} \quad \sqrt{\cos \theta} > 0.7 \quad \text{so that} \quad {\lambda \over \sqrt{\cos \theta}} < 4.$$
Then, if $n \geq 4m$, we let 
$$\delta={\omega \over m-1}, \quad \tau=\sqrt{\cos \theta}$$
and the proof proceeds as above.
{\hfill \hfill \hfill} \qed
\remark{(3.6) Remark} As follows from (3.5.1), we have $P_{\Omega}(Z) \ne 0$ for any symmetric complex matrix $Z=\left(z_{ij}\right)$ satisfying the conditions of Theorem 1.5. The algorithm of Section 2 extends in a straightforward way to computing 
$P_{\Omega}(W)$ for every matrix $W=\left(w_{ij}\right)$ satisfying (1.2.1). This allows us to compute in $n^{O(\ln m)}$ time an
$m$-subset $S$ which satisfies (1.1.2), by conditioning successively on $i \in S$ for $i=1, \ldots, n$.

\endremark

\head Acknowledgments \endhead

I am grateful to Alex Samorodnitsky and Benny Sudakov for several useful remarks.


\Refs
\widestnumber\key{AAAA}

\ref\key{Az88}
\by R. Azencott
\paper Simulated annealing
\paperinfo S\'eminaire Bourbaki, Vol. 1987/88
\jour Ast\'erisque 
\pages No. 161--162,  223--237 
\yr 1988
\endref

\ref\key{BG08}
\by  A. Bandyopadhyay and D. Gamarnik
\paper Counting without sampling: asymptotics of the log-partition function for certain statistical physics models
\jour Random Structures $\&$ Algorithms 
\vol 33 
\yr 2008
\pages no. 4, 452--479
\endref

\ref\key{B+07}
\by M. Bayati, D. Gamarnik, D. Katz, C. Nair and P.  Tetali
\paper Simple deterministic approximation algorithms for counting matchings
\inbook STOC'07 -- Proceedings of the 39th Annual ACM Symposium on Theory of Computing
\pages 122--127
\publ ACM
\publaddr New York
\yr 2007
\endref 

\ref\key{Bh12}
\by A. Bhaskara
\book Finding Dense Structures in Graphs and Matrices
\bookinfo Ph.D. dissertation, Princeton University, available at \hfill \hfill \hfill \break {\tt http://www.cs.princeton.edu/$\sim$bhaskara/thesis.pdf}
\yr 2012
\endref

\ref\key{B+10}
\by A. Bhaskara, M. Charikar, E. Chlamtac, U. Feige and A. Vijayaraghavan
\paper Detecting high log-densities -- an $O(n^{1/4})$ approximation for densest $k$-subgraph
\inbook STOC'10--Proceedings of the 2010 ACM International Symposium on Theory of Computing
\pages 201--210
\publ ACM
\publaddr New York
\yr 2010
\endref


\ref\key{BG05}
\by A. Bulatov and M. Grohe
\paper The complexity of partition functions
\jour Theoretical Computer Science
\vol 348 
\yr 2005
\pages no. 2--3, 148--186
\endref

\ref\key{F+01}
\by U. Feige, G. Kortsarz and D. Peleg
\paper The dense k-subgraph problem
\jour Algorithmica 
\vol 29 
\yr 2001
\pages no. 3, 410--421
\endref

\ref\key{FK99}
\by A. Frieze and R. Kannan
\paper Quick approximation to matrices and applications
\jour  Combinatorica 
\vol 19 
\yr 1999
\pages no. 2, 175--220
\endref

\ref\key{G+98}
\by O. Goldreich, S. Goldwasser and D. Ron
\paper Property testing and its connection to learning and approximation
\jour  Journal of the ACM 
\vol 45 
\yr 1998
\pages no. 4, 653--750
\endref

\ref\key{H\aa 99}
\by J. H\aa stad
\paper Clique is hard to approximate within $n^{1-\epsilon}$
\jour Acta Mathematica 
\vol 182 
\yr 1999
\pages  no. 1, 105--142
\endref

\ref\key{JS93}
\by M. Jerrum and A. Sinclair
\paper Polynomial-time approximation algorithms for the Ising model
\jour SIAM Journal on Computing 
\vol 22 
\yr 1993
\pages no. 5, 1087--1116
\endref 

\ref\key{LY52}
\by T.D. Lee and C.N. Yang
\paper Statistical theory of equations of state and phase transitions. II. Lattice gas and Ising model
\jour Physical Review (2) 
\vol 87
\yr 1952
\pages 410--419
\endref 

\ref\key{LP09}
\by L. Lov\'asz and M.D. Plummer
\book Matching Theory
\bookinfo Corrected reprint of the 1986 original
\publ AMS Chelsea Publishing
\publaddr Providence, RI
\yr 2009
\endref

\ref\key{SS05}
\by A.D. Scott and A.D. Sokal
\paper The repulsive lattice gas, the independent-set polynomial, and the Lov\'asz local lemma
\jour Journal of Statistical Physics
\vol 118 
\yr 2005
\pages no. 5-6, 1151--1261
\endref

\ref\key{Zu07}
\by D. Zuckerman
\paper Linear degree extractors and the inapproximability of max clique and chromatic number
\jour Theory of Computing
\vol  3 
\yr 2007
\pages 103--128
\endref



 
\endRefs
\enddocument

\end